\documentclass[12pt]{article}
\usepackage{amssymb}
\usepackage{graphicx}
\usepackage{amsfonts}
\usepackage{latexsym}
\usepackage{amsthm}
\usepackage{amsmath}

\usepackage{amssymb, amscd}

\usepackage{amsmath,amsfonts,graphicx,subfigure,url,amscd}

\usepackage{url}
\usepackage[T1]{fontenc}
\usepackage{color}
\usepackage{tabu}

\usepackage{bbm}

\numberwithin{equation}{section}

\newtheorem{conj}{Conjecture}[section]
\newtheorem{thm}[conj]{Theorem}
\newtheorem{rem}[conj]{Remark}

\newtheorem{defin}[conj]{Definition}

\newtheorem{cor}[conj]{Corollary}
\newtheorem{lema}[conj]{Lemma}

\begin{document}
\date{January 25, 2020}
 \title{{Generalized chessboard complexes \\ and discrete Morse theory\footnote{This research was supported by the Grants 174020 and
174034 of the Ministry of Education,
Science and Technological Development of the Republic of Serbia. The authors acknowledge the hospitality of Mathematical Institute in Oberwolfach, where this paper was completed.}}}

\author{{Du\v{s}ko Joji\'{c}} \\
 {\small Faculty of Science}\\[-2mm]
 {\small University of Banja Luka}
\and Gaiane Panina\\
{\small St. Petersburg State University} \\[-2mm]
{\small St. Petersburg Department of}\\[-2mm] {\small Steklov Mathematical Institute}
\and Sini\v{s}a T. Vre\'{c}ica\\
{\small Faculty of Mathematics}\\[-2mm] {\small University of Belgrade}
\and Rade T. \v{Z}ivaljevi\'{c}\\
{\small Mathematical Institute}\\[-2mm] {\small SASA,
  Belgrade}\\[-2mm]}

\maketitle
\begin{abstract}\noindent
Chessboard complexes and their generalizations, as objects, and Discrete Morse theory, as a tool, are presented as a unifying theme linking different areas of geometry, topology, algebra and combinatorics. Edmonds and Fulkerson bottleneck (minmax)  theorem is proved  and interpreted as a result about a critical point of a discrete Morse function on the Bier sphere $Bier(K)$ of an associated simplicial complex $K$. We illustrate the use of ``standard discrete Morse functions'' on generalized chessboard complexes by proving a connectivity result for chessboard complexes with multiplicities. Applications include new Tverberg-Van Kampen-Flores type results for $j$-wise disjoint partitions of a simplex.
\end{abstract}

\section{Introduction}
\label{sec:intro}

Chessboard complexes and their relatives have been
for decades an important theme of topological combinatorics. They have found numerous and
often unexpected applications in group theory, representation
theory, commutative algebra, Lie theory, discrete and computational geometry,
algebraic topology, and geometric and topological combinatorics, see \cite{Ata04}, \cite{AuFie07}, \cite{Au10} \cite{BLVZ}, \cite{FrHa98}, \cite{Ga79}, \cite{Jo07}, \cite{krw}, \cite{ShaWa04}, \cite{VZ94}, \cite{VZ07}, \cite{Wa03}, \cite{Zie94}, \cite{ZV92}.

The books \cite{Jo-book} and \cite{M}, as well as the review papers \cite{Wa03} and
\cite{Z17}, cover selected topics of the theory of chessboard
complexes and contain a more complete list of related
publications.

\medskip
Chessboard complexes and their generalizations are some of the
most studied {\em graph complexes} \cite{Jo-book}. From this point of view chessboard complexes can be interpreted as relatives of L. Lov\'asz $Hom$-complexes  \cite{Ko}, matching complexes,
clique complexes, and many other important classes of simplicial complexes.

\medskip
More recently new classes of generalized chessboard complexes have emerged and new methods, based on novel {\em shelling techniques}  and ideas from Forman's {\em discrete Morse theory}, were introduced. Examples include multiple and symmetric multiple chessboard complexes \cite{jvz, jvz2}, Bier complexes \cite{jnpz}, and deleted joins of collectively unavoidable complexes, see \cite{jnpz} and \cite{jpz1, jpz2}. Among applications are the resolution of the balanced case of the ``admissible/prescribable partitions conjecture'' \cite{jvz2}, general Van Kampen-Flores type theorem for balanced, collectively unavoidable complexes \cite{jpz1}, and
``balanced splitting necklace theorem'' \cite{jpz2}.

\medskip
This paper is both a leisurely introduction and an invitation to this part of topological combinatorics, and a succinct overview of some of the ideas of discrete Morse theory, combinatorics and equivariant topology, used in our earlier papers.

\medskip
New results are in Sections~\ref{sec:Ed-Fu}, \ref{sec:dusko} and \ref{sec:sinisa}.  They include an alternative treatment of Edmonds and Fulkerson bottleneck (minmax) theorem (Section~\ref{sec:Ed-Fu}) and the construction of ``standard discrete Morse functions'' on generalized chessboard complexes with multiplicities (Section~\ref{sec:dusko}). This leads to a frequently optimal connectivity result for generalized chessboard complexes with multiplicities (Theorem~\ref{T:connect} in Section~\ref{sec:dusko}), which is  used in Section~\ref{sec:sinisa} for the proof of new Tverberg-Van Kampen-Flores type results for $j$-wise disjoint partitions of a simplex.

\section{Chessboard complexes}
\label{sec:chessboard}

Chessboard complexes naturally arise in the study of the geometry of {\em admissible rook configurations} on a general $(m\times n)$-chessboard. An admissible configuration is any non-taking placement of rooks, i.e., a placement which does not allow any two of them to be in the same row or in the same column. The collection of all these placements forms a
simplicial complex which is called the chessboard complex and denoted by $\Delta_{m,n}$.

\medskip
More formally, the set of vertices of $\Delta_{m,n}$ is ${\rm Vert}(\Delta_{m,n}) = [m]\times [n]$ and $S\subseteq [m]\times [n]$ is a simplex of  $\Delta_{m,n}$ if and only if for each two distinct elements $(i_1, j_1), (i_2, j_2) \in S$ neither $i_1 = i_2$ nor $j_1 = j_2$.

\subsection{An example}

\begin{figure}[htb]
\centering %\vspace{-2cm}
 \includegraphics[scale=.80]{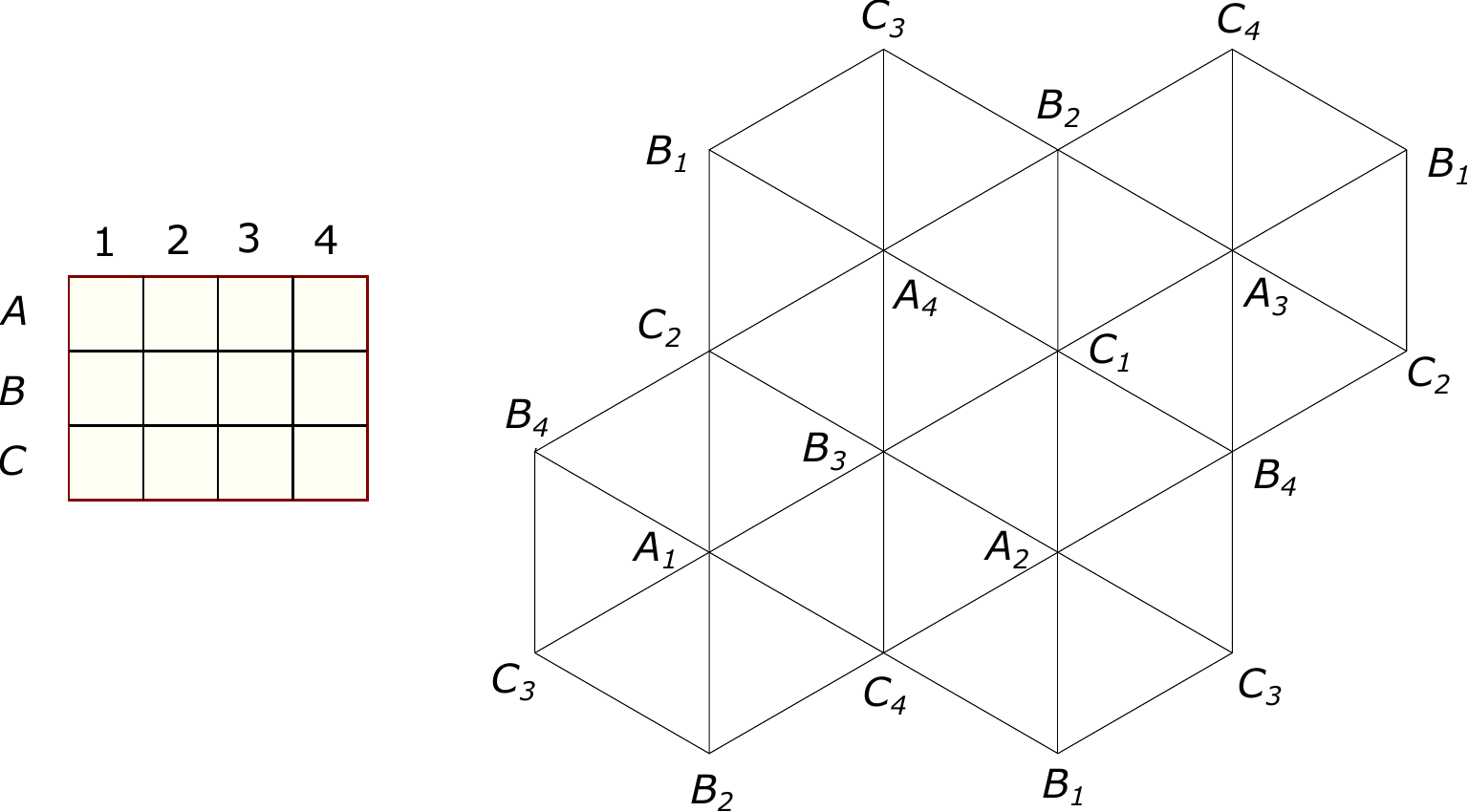}
\caption{Chessboard complex $\Delta_{4,3}$ } \label{fig:delta_4_3}
\end{figure}

Let us take a closer look at one of the simplest chessboard complexes, the complex $\Delta_{4,3}$, based
on the $4\times 3$ chessboard (see Figure~\ref{fig:delta_4_3}).

The $f$-vector of $\Delta_{4,3}$ is $f(\Delta_{4,3}) = (12, 36, 24)$ so its Euler characteristics is $\chi(\Delta_{4,3})=0$. Moreover, the geometric realization of $\Delta_{4,3}$ is an orientable $2$-dimensional manifold.

Indeed, the link of each vertex is isomorphic to $\Delta_{3,2}$ (= hexagonal triangulation of the circle $S^1$) while the link of each edge is the circle $S^0$. Each $2$-dimensional simplex $\sigma = \{A_i,B_j,C_k\}$ is uniquely completed to a permutation $\pi = (i,j,k,l)$ of the set $[4] = \{1,2,3,4\}$ and ${\rm Sign}(\sigma) := {\rm Sign}(\pi)$ defines an orientation on $\Delta_{4,3}$.

\medskip
From here we immediately conclude that $\Delta_{4,3}$ is a triangulation of the $2$-dimensional torus $T^2$. The universal covering
of $\Delta_{4,3}$ is identified as the honeycomb tiling of the plane and the corresponding fundamental domain is exhibited in
Figure~\ref{fig:delta_4_3}. From here we can easily read off the generators of the group $H_1(\Delta_{4,3}; \mathbb{Z}) \cong
\mathbb{Z}^2$ as the geodesic edge-paths connecting the three copies of vertex $C_3$, shown in Figure~\ref{fig:delta_4_3}.

\subsection{Graph complexes}

Let $G$ be a finite graph with vertex set $V = V_G$ and edge set $E= E_G$. A graph complex
on $G$ is an abstract simplicial complex consisting of subsets of $E$. We usually interpret such a complex as a family of subgraphs of $G$. The study of graph complexes,  with the  emphasis on their homology,
homotopy type, connectivity degree, Cohen-Macaulayness, etc., has been an active area of study in topological combinatorics, see
\cite{Jo-book}.

\medskip
The chessboard complex $\Delta_{m,n}$ can be interpreted as a graph complex of the complete bipartite graph $K_{m,n}$, where
the simplices $S\subset [m]\times [n]$ are interpreted as ``matchings'' in $K_{m,n}$. Recall that $\Gamma \subseteq E_G$ is
a matching on the graph $G$ if each $v\in V_G$ is incident to at most one edge in $\Gamma$.

\medskip
All ``generalized chessboard complexes'', introduced in Section~\ref{sec:gen-chessboard}, can be also described as graph
complexes of the  graph $K_{m,n}$.

\subsection{Chessboard complexes as Tits coset complexes}

Perhaps the first appearance of chessboard complexes was in the thesis of Garst \cite{Ga79}, as {\em Tits coset complexes}.
Recall that a Tits coset complex $\Delta(G; G_1,\dots, G_n)$, associated to a group $G$ and a family $\{G_1,\dots, G_n\}$ of its subgroups is the nerve $Nerve(\mathcal{F})$ of the associated family of cosets $\mathcal{F} = \{gG_i \mid g\in G, i\in [n]\}$.
More explicitly vertices of $\Delta(G; G_1,\dots, G_n)$ are cosets $gG_i$ and a colection $S = \{g_jG_i\}_{(i,j)\in I}$, for some
$I\subseteq [n]  \times G$, is a simplex of $\Delta(G; G_1,\dots, G_m)$ if and only if
\[
     \bigcap_{(i,j)\in I}  g_jG_i \neq \emptyset \, .
\]
If $G = S_m$ is the symmetric group and $G_i := \{\pi\in S_m \mid \pi(i) = i\}$   for $i=1,\dots, n$, the associated Tits coset complex is the chessboard complex $\Delta_{m,n}$.

\subsection{Chessboard complexes in discrete geometry}

Chessboard complexes made their first appearance in discrete geometry in \cite{ZV92}, in the context of the so called {\em colored Tverberg problem}.

For illustration, an instance of the {\em type~B} colored Tverberg theorem \cite{VZ94, Z17} claims that
for each collection $C\subset \mathbb{R}^3$ of fifteen points in the $3$-space, evenly colored by three colors, there exist three
vertex disjoint triangles $\Delta_1, \Delta_2, \Delta_3$, formed by the points of different color, such that $\Delta_1\cap\Delta_2\cap\Delta_3\neq\emptyset$.

\medskip
A general form of this result was deduced in \cite{VZ94} from a Borsuk-Ulam type result claiming that each $\mathbb{Z}_r$-equivariant map
\begin{equation}\label{eqn:nasa}
(\Delta_{r,2r-1})^{\ast
(k+1)}\stackrel{\mathbb{Z}_r}{\longrightarrow} W_r^{\oplus (d+1)}
\end{equation}
 must have a zero if $r\leq d/(d-k)$ (this is a necessary condition), $r$ is a prime power, $\Delta_{r,2r-1}$ is a chessboard complex,
 and $W_r = \{x\in \mathbb{R}^r \mid x_1+\dots + x_r =0\}$.

 The reader is referred to \cite{Z17} for an overview of these and more recent results, as well as for a more complete list of references.

\section{Generalized chessboard complexes}
\label{sec:gen-chessboard}

Motivated primarily by applications to problems in discrete geometry, especially the problems of Tverberg and Van Kampen-Flores type, more general chessboard complexes were introduced and studied. Closely related complexes previously emerged in algebraic combinatorics \cite{krw, Wa03}.

These complexes are also referred to as  {\em generalized chessboard complexes}, since the set of vertices remains the $(m\times n)$-chessboard $[m]\times [n]$, but the criterion for $S\subseteq [m]\times [n]$ to be a simplex (``admissible rook placement'') may be quite different and vary from problem to problem.

\medskip
The following definition includes most if not all of the currently studied examples and provides a natural ecological
niche for all these complexes and their relatives.

\begin{defin}\label{def:general-chessboard}
Suppose that $\mathcal{K} = \{K_i\}_{i=1}^n$ and $\mathcal{L} =
\{L_j\}_{j=1}^m$ are two collections of simplicial complexes where
${\rm Vert}(K_i) = [m]$ for each $i\in [n]$ and ${\rm Vert}(L_j) =
[n]$ for each $j\in [m]$. Define,
\begin{equation}\label{eqn:2-collections}
\Delta_{m,n}^{\mathcal{K}, \mathcal{L}} =
\Delta_{m,n}(\mathcal{K}, \mathcal{L})
\end{equation}
as the complex of all subsets (rook-placements) $A\subset
[m]\times [n]$ such that $\{i\in [m]\mid (i,j)\in A \}\in K_j$ for
each $j\in [n]$ and $\{j\in [n]\mid (i,j)\in A \}\in L_i$ for each
$i\in [m]$.
\end{defin}

Definition~\ref{def:general-chessboard} can be specialized in many
ways. Again, we focus on the special cases motivated by intended
applications to the generalized Tverberg problem.

\begin{defin}\label{def:special-1}
Suppose that $\mathbf{k}=(k_i)_{i=1}^n$ and
$\mathbf{l}=(l_j)_{j=1}^m$ are two sequences of non-negative
integers. Then the complex,
\begin{equation}\label{eqn:defin-hv}
\Delta_{m,n}^{\mathbf{k}, \mathbf{l}} =
\Delta_{m,n}^{k_1,\dots,k_n;l_1,\dots,l_m}
\end{equation}
arises as the complex of all rook-placements $A\subset [m]\times
[n]$ such that at most $k_i$ rooks are allowed to be in the $i$-th
row (for $i=1,\dots,n$), and at most $l_j$ rooks are allowed to be
in the $j$-th column (for $j=1,\dots,m$).
\end{defin}
\begin{rem} {\rm
The complexes $\Delta_{m,n}^{\mathbf{k}, \mathbf{l}} = \Delta_{m,n}^{k_1,\dots,k_n;l_1,\dots,l_m}$ are sometimes referred to
as the {\em chessboard complexes with multiplicities} or {\em multiple chessboard complexes}. Closely related are ``bounded degree
graph complexes'', studied in \cite{krw} and \cite{Wa03}.  }
\end{rem}

When $k_1=\dots =k_n=p$ and $l_1=\dots =l_m=q$, we obtain the
complex $\Delta_{m,n}^{p,q}$. For the reasons which will become
clear in the following section of the paper, in our earlier papers \cite{jvz, jvz2} we focused to
 the case $l_1=\dots =l_m=1$, i.e. to the complexes,
\begin{equation}\label{eqn:oznaka}
\Delta_{m,n}^{k_1,\dots,k_n;\mathbf{1}} :=
\Delta_{m,n}^{k_1,\dots,k_n;1,\dots,1} \, .
\end{equation}
In Section~\ref{sec:dusko} of this paper we fill this ``gap'' and return to the case of general chessboard complexes with multiplicities.

\subsection{$n$-fold $j$-wise deleted join}

Joins and deleted joins of simplicial complexes, as well as their generalizations, have found numerous applications in
topological combinatorics, see \cite[Section 6.3]{M} for motivation and an introduction.

\medskip
For a simplicial complex $K$, the $n$-fold $j$-wise deleted join of $K$ is
\begin{equation}\label{eqn:j-wise}
K^{\ast n}_{\Delta(j)} := \{A_1\uplus A_2\uplus\dots\uplus A_n \in K^\ast \mid (A_1,A_2,\dots, A_n) \mbox{ {is $j$-wise disjoint}}\}
\end{equation}
where an $n$-tuple $(A_1,A_2,\dots, A_n)$ is $j$-wise disjoint if every sub-collection $\{A_{k_i}\}_{i=1}^j$, where
$k_1<k_2<\dots <k_j$, has an empty intersection.

\medskip
It immediately follows that $K^{\ast n}_{\Delta(j)} \subseteq K^{\ast n}_{\Delta(j+1)}$ and that  $K^{\ast n}_{\Delta(n+1)} = K^{\ast n}$ and  $K^{\ast n}_{\Delta(2)} = K^{\ast n}_{\Delta}$ are
respectively the $n$-fold join and the $n$-fold deleted join of the complex $K$.

\medskip
A simple but very useful property of these operations is that they commute
\[
     (K^{\ast n}_{\Delta(j)})^{\ast m}_{\Delta(k)} \cong (K^{\ast m}_{\Delta(k)})^{\ast n}_{\Delta(j)} \, .
\]
For example if $K = pt$ is a one-point simplicial complex we obtain the isomorphsim
\[
    \Delta_{m,n} =  ((pt)^{\ast m}_{\Delta})^{\ast n}_\Delta \cong  ((pt)^{\ast n}_{\Delta})^{\ast m}_\Delta = \Delta_{n,m} \, .
\]
A single complex $K$ in equation (\ref{eqn:j-wise}) can be replaced by a collection $\mathcal{K} = \{K_j\}_{j=1}^n$ of complexes
$K_j\subseteq 2^{[m]}$ which leads to the definition of the $j$-wise deleted join of $\mathcal{K}$,
\begin{equation*}
\mathcal{K}_{\Delta(j)}  := \{A_1\uplus \dots\uplus A_n \in K_1\ast\dots\ast K_n \mid (A_1,\dots, A_n) \mbox{ {is $j$-wise disjoint}}\} \, .
\end{equation*}

All simplicial complexes described in this section are generalized chessboard complexes in the sense of Definition~\ref{def:general-chessboard}. For example if $K\subseteq 2^{[m]}$ then its $n$-fold $j$-wise deleted join is the complex
\[
 K^n_{\Delta(j)} \cong  \Delta_{m,n}^{\mathcal{K}, \mathcal{L}}
\]
where $K_1=\dots = K_n$ and $L_1=\dots = L_m = {{[m]}\choose{\leq j-1}}$ is the collection of all subsets of $[m]$ of cardinality strictly less than $j$.

\subsection{Bier spheres as generalized chessboard complexes}

Let $K \varsubsetneq 2^{[m]}$ be a simplicial complex on the ground set
$[m]$ (meaning that we allow $\{j\}\notin K$ for some $j\in[m]$).
The Alexander dual of $K$ is the simplicial complex $K^\circ$ that
consists of the complements of all nonsimplices of $K$
\[
K^\circ := \{A^c \mid A\notin K \} \, .
\]
By definition the {\em ``Bier sphere''} is the deleted join $Bier(K) := K \ast_\Delta K^\circ$. (A face $A_1\uplus A_2\in Bier(K)$ is often denoted as a triple $(A_1, A_2; B)$ where $B:= [m]\setminus (A_1\cup A_2)$.)

It turns out that $Bier(K)$ is indeed a triangulation of an $(m-2)$-dimensional sphere \cite{Bier},
 see \cite{M} and \cite{long-1} for different,  very short and elegant proofs.

\medskip
The Bier sphere $Bier(K)$ is also a generalized chessboard complex where $K_1= K, K_2 = K^\circ$ and $L_1=\dots = L_m = \{\emptyset, \{1\}, \{2\}\}\subset 2^{[2]}$.

\medskip
{\em Alexander $r$-tuples} $\mathcal{K}
= \{K_i\}_{i=1}^r$ of simplicial complexes were introduced in \cite{jnpz} as a
generalization of pairs $(K, K^\circ)$  of Alexander dual complexes. The associated
 {\em  generalized Bier complexes}, defined as the $r$-fold deleted joins
$\mathcal{K}^{\ast r}_\Delta$ of Alexander $r$-tuples are also generalized chessboard complexes in the sense of Definition~\ref{def:general-chessboard}.

\section{Discrete Morse theory}
\label{sec:DMT}

A discrete Morse function on a simplicial complex $K\subseteq 2^V$ is, by definition, an acyclic
matching on the Hasse diagram of the partially ordered set $(K, \subseteq)$. Here is a brief reminder of the basic facts and definitions of discrete Morse theory.

\medskip

Let $K$ be a simplicial complex. Its $p$-dimensional simplices
($p$-simplices for short) are denoted by $\alpha^p, \alpha^p_i, \beta^p,
\sigma^p$, etc. A \textit{discrete vector field} is a set of
pairs $D = \{\dots, (\alpha^p,\beta^{p+1}), \dots\}$ (called a matching) such
that:
\begin{enumerate}
 \setlength\itemsep{-1.5mm}
    \item[(a)]  each simplex of the complex participates in at most one
    pair;
    \item[(b)]  in each pair $(\alpha^p,\beta^{p+1})\in D$, the simplex $\alpha^p$ is a facet of $\beta^{p+1}$;
    \item[(c)]  the empty set $\emptyset\in K$ is not
matched, i.e.\ if $(\alpha^p,\beta^{p+1})\in D$ then $p\geq 0$.
\end{enumerate}
\noindent The pair $(\alpha^p, \beta^{p+1})$ can be informally
thought of as a vector in the vector field $D$. For this reason it
is occasionally denoted by $\alpha^p \rightarrow \beta^{p+1}$ or $\alpha^p \nearrow \beta^{p+1}$  (and
in this case $\alpha^p$ and $\beta^{p+1}$ are informally referred to as  the {\em beginning} and {\em the end}\/ of
the arrow $\alpha^p \rightarrow \beta^{p+1}$).

\smallskip
Given a discrete vector field $D$, a \textit{gradient path} in $D$
is a sequence of simplices (a zig-zag path)
$$\alpha_0^p \nearrow \beta_0^{p+1} \searrow \alpha_1^p \nearrow \beta_1^{p+1} \searrow \alpha_2^p \nearrow \beta_2^{p+1} \searrow \,  \cdots \,   \searrow\alpha_m^p \nearrow \beta_m^{p+1} \searrow \alpha_{m+1}^p$$
satisfying the following conditions:
\begin{enumerate}\setlength\itemsep{-1.5mm}
\item  $\big(\alpha_i^p,\ \beta_i^{p+1}\big)$ is a pair in $D$ for each $i$;
\item for each $i = 0,\dots, m$ the simplex $\alpha_{i+1}^p$
 is a facet of $\beta_i^{p+1}$;
 \item for each $i = 0,\dots, m-1,
 \, \alpha_i\neq \alpha_{i+1}$.
\end{enumerate}

A path is \textit{closed} if $\alpha_{m+1}^p=\alpha_{0}^p$. A
\textit{discrete Morse function } (DMF for short) is a discrete
vector field without closed paths.

\medskip
Assuming that a discrete Morse function is fixed, the {\em
critical simplices} are those simplices of the complex that are
not matched. The Morse inequality \cite{Forman2} implies that critical
simplices cannot be completely avoided.

A discrete Morse function $D$ is \textit{perfect}
if the number of critical $k$-simplices  equals the $k$-th
Betty number of the complex. It follows that $D$ is a perfect Morse function if and only if the number of all critical simplices equals the sum of all Betty numbers of $K$.

\medskip
A central idea of discrete Morse theory, as summarized in
the following theorem of R.~Forman, is to contract all matched
pairs of simplices and to reduce the simplicial complex $K$ to a
cell complex (where critical simplices correspond to the cells).

\begin{thm}\cite{Forman2} \label{ThmWedge} Assume that a
discrete Morse function on a simplicial complex $K$ has a single
zero-dimensional critical simplex $\sigma^0$ and that all other
critical simplices have the same dimension $N>1$. Then $K$ is
homotopy equivalent to a wedge of $N$-dimensional spheres.

More generally, if all critical simplices, aside from $\sigma^0$, have
dimension $\geq N$, then the complex $K$ is $(N-1)$-connected.  \qed
\end{thm}

\subsection{Discrete vector fields on Bier spheres}
\label{sec:DMF-Bier}

It is known that all Bier spheres are shellable, see
\cite{BjornerZiegler} and \cite{cukic}. A method of Chari \cite{Chari} can be used
to turn this shelling into a perfect discrete Morse function (DMF). The
construction of our {\em perfect DMF} on a Bier sphere
 essentially follows this path, see \cite{jnpz} for more details. For the reader's convenience here we reproduce this
 construction since it will be needed in Section~\ref{sec:Ed-Fu}.

\bigskip
{\bf A perfect DMF on $Bier(K)$}

\medskip
We construct a discrete vector field $D_1$ on the Bier sphere
$Bier(K)$ in two steps:
\begin{enumerate}\setlength\itemsep{-1.5mm}
    \item[(1)] We match  the simplices
$$\alpha=(A_1,A_2;B\cup i)  \hbox{ and } \beta=(A_1,A_2\cup i;B )$$ iff the following holds:
\begin{enumerate}\setlength\itemsep{-1.5mm}
    \item[(i)] $i<B,\ i<A_2$  \newline(that is, $i$ is smaller than all the entries of $B$ and $A_2$).
    \item[(ii)] $A_2\cup i \in K^{\circ}$.
\end{enumerate}
\end{enumerate}

\medskip\noindent
Before we pass to step 2, let us observe that the non-matched
simplices are labelled  by $(A_1,A_2;B\cup i)$ such that $A_2\in
K^{\circ}$, but $A_2\cup i \notin K^{\circ}$. As a consequence,
for non-matched simplices $A_1\cup B\in K$.

\smallskip
\begin{enumerate}\setlength\itemsep{-1.5mm}
    \item[(2)] In the second step we match together the simplices
$$\alpha=(A_1,A_2;B\cup j)  \hbox{ and } \beta=(A_1\cup j,A_2;B)$$ iff the following holds:
\begin{enumerate}\setlength\itemsep{-1.5mm}
\item None of the simplices $\alpha$ and $\beta$ is matched in the
first step.
    \item $j>B,\ j>A_1$.
    \item $A_1\cup j \in K$.
\end{enumerate}
\end{enumerate}
\noindent Observe that the condition (c)  always holds (provided
that the condition (a) is satisfied), except for
    the case  $B=\emptyset$.

\begin{lema}{\rm (see \cite[Lemma 6.1]{jnpz})} The discrete vector field $D_1$ is a discrete Morse function on the Bier sphere $Bier(K)$.
\end{lema}

\textit{Proof.} Since $D_1$ is (by construction) a discrete vector
field,  it remains to check that there are no closed gradient
paths. Observe that in each pair of simplices in the discrete
vector field $D_1$ there is exactly one \textit{migrating
element}. More precisely, in the case (1) the element $i$
{migrates} to $A_2$, and in  the case (2) the element $j$
{migrates} to $A_1$.

The lemma follows from the observation that (along a gradient
path) the values of the migrating element that move to $A_2$
strictly decreases. Similarly, the values of migrating elements
that move to $A_1$ can only increase.
This is certified through the following  simple case analysis:
(1) After a first step pairing  comes a splitting of $A_2$. Then the gradient path terminates.
(2) After a first step pairing (with migrating element $i$)  comes a splitting of $A_1$. The gradient path proceeds only if the  splitted element is smaller than $i$. (2) After a second step pairing  comes a splitting of $A_1$. Then the gradient path terminates.
(2) After a second step pairing (with migrating element $i$)  comes a splitting of $A_2$. The gradient path proceeds only if the  splitted element is bigger than $i$.

\medskip
Let us illustrate this observation by an example which captures the above case analysis. Assume we have a fragment of a gradient path that contains two matchings of type 1.
We have:

$$(A_1\cup k ,A_2;B\cup i ) \rightarrow (A_1\cup k ,A_2\cup i;B)\rightarrow $$$$(A_1,A_2\cup i;B\cup k ) \rightarrow (A_1,A_2\cup k \cup i;B)$$
The migrating elements here are $i$ and $k$. The definition of the
matching $D_1$ implies $k < i$. Otherwise $(A_1,A_2\cup i;B\cup k
)$ is matched with $(A_1,A_2;B\cup k\cup i )$, and the path would
terminate after its second term. \qed

\medskip

It is not difficult to see that there are precisely two critical
simplices in $D_1$:

\begin{enumerate}\setlength\itemsep{-1.5mm}
    \item An $(n-2)$-dimensional simplex,  $$(A_1,A_2; i )$$
where $A_1<i<A_2$,  (this condition describes this simplex
uniquely, in light of the fact that $A_1\in K$ and $A_2\in
K^{\circ}$),
    \item and the $0$-dimensional simplex, $$(\emptyset, \{1\} ;\{2,3,4,...,n\}).$$
\end{enumerate}
\noindent (Here we make a simplifying assumption that $\{1\}\in
K^\circ$, which can be always achieved by a re-enumeration, except
in the trivial case $K^\circ =\{\emptyset\}$.)

\subsection{Discrete vector fields on generalized chessboard complexes}

The construction of the discrete Morse function on the Bier sphere $Bier(K)$ illustrates the fruitful idea which
  can be extended and further developed to cover the case of other generalized chessboard complexes.

\medskip Examples of this construction can be found in \cite{jnpz} and \cite{jpz1}, see also Section~\ref{sec:dusko} for a construction of such a discrete Morse function on the multiple chessboard complex  $\Delta_{m,n}^{k_1,\dots,k_n;l_1,\dots,l_m}$.

\medskip All these constructions of DMF share the same basic idea,  for this reason we sometimes refer to them as {\em standard DMF on generalized chessboard complexes}. Note that the proofs that they indeed form an acyclic matching may vary from example to example and use some special properties of the class under investigation.

\section{Edmonds-Fulkerson bottleneck extrema}
\label{sec:Ed-Fu}

  In this section we connect, via discrete Morse theory, the combinatorial topology
  of Bier spheres with Edmonds-Fulkreson theorem on bottleneck extrema of pairs of dual clutters.
  We will show that there is much more than meets the eye in the standard concise treatment of this classical result of
  combinatorial optimization.

\begin{figure}[htb]
\centering %\vspace{-2cm}
 \includegraphics[scale=.80]{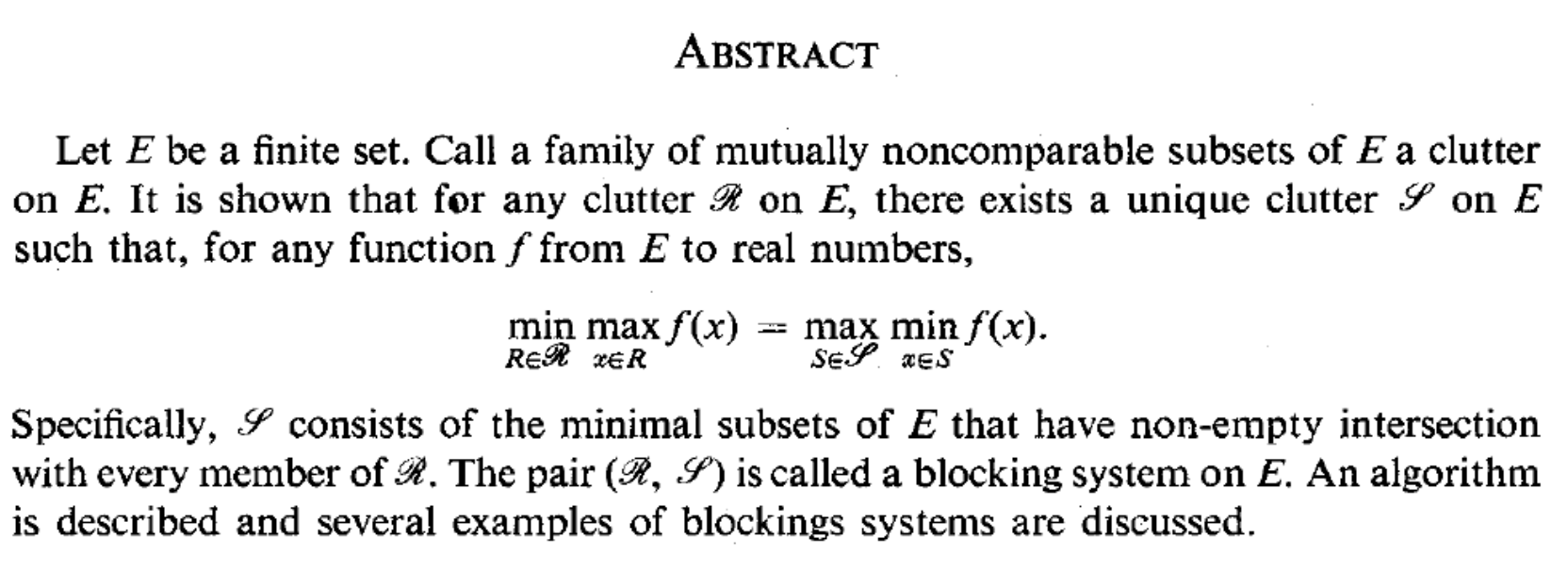}
\caption{Edmonds-Fulkerson bottleneck theorem}\label{fig:bottleneck}
\end{figure}

\medskip
Figure~\ref{fig:bottleneck} shows the abstract of the published version of \cite{Ed-Fu}, which originally appeared as a RAND-corporation preprint AD 664879 in January of 1966.

This is a purely combinatorial result which is often referred to as the Edmonds-Fulkreson bottleneck lemma (theorem). Minmax theorems are ubiquitous in mathematics, notably in geometry, polyhedral combinatorics, critical point theory, game theory and other areas. One of early examples is the minimax theorem of John von Neumann (first proven and published in 1928) which gives conditions on a function $f : C\times D \rightarrow \mathbb{R}$, defined on
the product of two closed, convex sets in $\mathbb{R}^n$, to satisfy the minmax equality,
\begin{equation}\label{eqn:minmax}
\min_{y\in D}\max_{x\in C} f(x,y) = \max_{x\in C}\min_{y\in D} f(x,y) \, .
\end{equation}
 It is interesting to compare the Edmonds-Fulkerson minmax theorem with their geometric counterparts. For example in a vicinity of a non-degenerate critical point a Morse function has the form
 $f(x,y) = -\vert x\vert^2 + \vert y \vert^2 = -x_1^2-\dots -x_p^2 + y_1^2+\dots + y_q^2$. Moreover, this function satisfies the concave/convex condition of  von Neumann's minmax theorem and the relation (\ref{eqn:minmax}) is valid.

\medskip
 There is a formal resemblance of these results, for example the $x$-sections (respectively $y$-sections) of the convex sets $C\times D$ in (\ref{eqn:minmax}) formally play the role of complementary clutters $\mathcal{R}$
 and $\mathcal{S}$ from the result of Edmonds and Fulkerson.   At first sight it appears to be naive and hard to expect a deeper connection between these results. Indeed, the clutter $\{C\times \{y\}\}_{y\in D}$ of $y$-sections is nowhere near to be
 the complementary clutter of the set $\{\{x\}\times D\}_{x\in C}$ of all $x$-sections, which is a consequence of the following
 lemma (see the property (3) on page 301 in \cite{Ed-Fu}).

 \begin{lema}\label{lema:complementary}
 The clutter $\mathcal{S}\subset 2^E$ is the complementary clutter of the clutter $\mathcal{R}\subset 2^E$, if and only if for each partition $E = E_0\uplus E_1$ of $E$ either an element of $\mathcal{R}$ is contained in $E_0$ or an element of $\mathcal{S}$ is contained in $E_1$, but not both.
 \end{lema}

 In the next section we show that there does exist a geometric interpretation of the Edmonds-Fulkerson bottleneck minmax equality, provided we are willing to replace the smooth by discrete Morse theory.

\subsection{Edmonds-Fulkerson minmax lemma revisited}

Here we use the results from Section~\ref{sec:DMF-Bier} to give a new proof and a new interpretation of Edmonds-Fulkerson minmax lemma. As before (Figure~\ref{fig:bottleneck}) the clutters $\mathcal{R}$ and $\mathcal{S}$ are both subfamilies of $2^E$.

\medskip

Let $\widehat{\mathcal{R}} := \{A\subseteq E \mid (\exists X\in \mathcal{R})\, X\subseteq A\}$ be the upper closure of the clutter $\mathcal{R}$ and let $K := 2^E\setminus \widehat{\mathcal{R}}$ be the complementary simplicial complex.

\begin{lema}
Let $K^\circ$ be the Alexander dual of the simplicial complex $K := 2^E\setminus \widehat{\mathcal{R}}$. Then
\[
        K^\circ = 2^E\setminus \widehat{\mathcal{S}}
\]
is the complementary simplicial complex of  the upper closure $\widehat{\mathcal{S}}$ of the clutter $\mathcal{S}$.
\end{lema}

\medskip\noindent
{\bf Proof:} This is an immediate consequence of Lemma~\ref{lema:complementary} since the pair of complexes $(K, K^\circ)$ is also characterized by the property that for each partition $E = E_0\uplus E_1$ precisely one of the following two relations
$E_0\in K, \, E_1\in K^\circ$ is satisfied. \hfill $\square$

\medskip Let $f : E \rightarrow \mathbb{R}$ be a real function. We may assume that $f$ is $1$-$1$. Moreover, we may replace $E$ by the set $[n]$ (where $n$ is the cardinality of $E$) and assume that $f = id : [n] \rightarrow [n]$ is the identity function.

\medskip
By construction and properties of the perfect DMF on the Bier sphere $Bier(K) = K\ast_\Delta K^\circ$, constructed in Section~\ref{sec:DMF-Bier}, there is a unique $(n-2)$-dimensional critical simplex $(A_1,A_2; i )$, characterized by the conditions
$A_1<i<A_2, \,  A_1\in K, \, A_2\in K^{\circ}$. Let us show that
\begin{equation}
  a:= \min_{I\in \mathcal{R}}\max_{x\in I} f(x)  = f(i) = \max_{J\in \mathcal{S}}\min_{x\in J} f(x) =: b \, .
\end{equation}
Indeed, $A_1\cup\{i\} \notin K$ implies $A_1\cup\{i\} \in \mathcal{R}$ and from $\max_{x\in A_1\cup\{i\}} f(x) = f(i)$ we
deduce the relation $a\leq f(i)$.

For the opposite inequality observe that if $I\in \mathcal{R}$ then $I\cap (A_2\cup \{i\}) \neq \emptyset$ (otherwise, since $A_2\cup \{i\}\in \mathcal{S}$, Lemma~\ref{lema:complementary} would be violated). Hence, $\max_{x\in I} f(x) \geq f(i)$ and $a\geq f(i)$.

The proof of the equality $b = f(i)$ is similar. \hfill $\square$

\begin{rem}{\rm One of the consequences is that the (algorithmic) complexity of determining the critical cell $(A_1, A_2; i)$ in the Bier sphere $Bier(K)$ is at least as big as the complexity of evaluating the maxmin (minmax) of a function on a family of sets
(clutter).  }
\end{rem}

 \section{Discrete Morse theory for chessboard complexes with multiplicities}
 \label{sec:dusko}

Suppose that $k_1,\ldots,k_n$ and $l_1,\ldots,l_m$ are two
sequences of non-negative integers. The generalized chessboard
complex $\Delta^{k_1,\ldots,k_n;l_1,\ldots,l_m}_{m,n}$ contains
all rooks placements on $[n]\times [m]$ table such that
at most $k_i$ rooks are in the $i$-th row and at most $l_j$ rooks
are in the $j$-th column. We use Forman's discrete Morse theory to
obtain a generalization of Theorem 3.2 from \cite{jvz}.
\begin{thm}\label{T:connect}
If $$l_1+l_2+\cdots+l_m\geqslant k_1+k_2+\cdots+k_n+n-1\ \ \ \ \ \ (\ast)$$ then
$\Delta^{k_1,\ldots,k_n;l_1,\ldots,l_m}_{m,n}$ is
$(k_1+k_2+\cdots+k_n-2)$-connected.
\end{thm}

\medskip\noindent
{\bf Proof:} A column (or a row) is called \textit{full} if it contains the maximal allowed number of rooks.
Otherwise, it is called \textit{free}.

We now define a Morse matching for
$\Delta=\Delta^{k_1,\ldots,k_n;l_1,\ldots,l_m}_{m,n}$. For a given face
$R$ we describe a face $R'$ that is paired with $R$, or we
recognize that $R$ is a critical face. Let us  do it stepwise.

\medskip
\textbf{Step 1.}

  Take the minimal $a_1$
 such that either  (1) there is a rook positioned at $(1,a_1)$,
 or (2) the $a_1$ column is free.

 In the first case (there is a rook at $(1,a_1)$), we match $R$ and
     $R'=R\setminus\{(1,a_1)\}$.

     This is always possible except for the unique exception, when $R$ contains exactly one rook at $(1,1)$.

     In the second case we match  $R$ and $R'=R\cup\{(1,a_1)\}$ provided that  $R'$ belongs to $\Delta$.
     The latter condition means that the first row in $R$ is not full.

     Clearly, after Step 1 the unmatched simplices are those with full first row, empty $(1, a_1)$,
      and a free column $a_1$.

\medskip
 \textbf{Step 2.} We match some of the simplices that are unpaired on the first step.

\begin{enumerate}\setlength\itemsep{-1.5mm}

\item If  there is  a rook at $(2,a_1)$,
           set $a_2:=a_1$ and match $R$ and $R'=R\setminus \{(2,a_2)\}$.
  \item If \begin{enumerate}\setlength\itemsep{-1.5mm}
             \item there is no  rook at $(2,a_1)$, and
             \item the number of rooks in column $a_{1}$ is smaller than $l_{a_1}-1$,
           \end{enumerate}

           set $a_2:=a_1$ and match $R$ and $R'=R\cup\{(2,a_2)\}$ provided that  $R'$ belongs to $\Delta$.
     The latter condition means that the second row in $R$ is not full.

     Introduce also $T(R):=2$.  Its meaning is  "the column $a_1=a_2$ has been used twice".

  \item If none of the above cases holds, set $a_2>a_1$ to be the minimal number such that
  either (1)   there is a rook positioned at $(2,a_2)$,
 or (2) the $a_2$ column is free.

The condition $(\ast)$  guarantees that $a_2$ is well-defined.

 If there is a rook at $(2,a_2)$, we match $R$ and
     $R'=R\setminus\{(2,a_2)\}$.

     Otherwise, we match  $R$ and $R'=R\cup\{(2,a_2)\}$ provided that  $R'$ belongs to $\Delta$.
     The latter condition means that the second row in $R$ is free.

     In this case we set $T(R):=1$, since the column $a_2$ has been used once.
\end{enumerate}

     Clearly, after Step 2 the unmatched simplices are those with full first and second rows, empty $(2, a_2)$,
      and  a free column $a_2$.

 \medskip

 We proceed in the same manner. During the first $k-1$ steps, some of the simplices become matched.
 Unmatched simplices have first $k-1$ rows full. They also have no rook at $(k-1,a_{k-1})$. Each unmatched simplex $R$ is associated a number $T(R)$.

\medskip
 This is how a generic step looks like:

 \textbf{Step k.}

 \begin{enumerate}\setlength\itemsep{-1.5mm}

 \item If  there is  a rook at $(k,a_{k-1})$,
          then  match $R$ and $R'=R\setminus \{(k,a_k)\}$.
  \item If \begin{enumerate}\setlength\itemsep{-1.5mm}
             \item there is no rook at  $(k,a_{k-1})$, and
             \item the number of rooks in column $a_{k-1}$ is smaller than $l_{a_{k-1}}-T(R)$,
           \end{enumerate}

           set $a_k:=a_{k-1}$ and match $R$ and $R'=R\cup\{(k,a_k)\}$ provided that  $R'$ belongs to $\Delta$.
     The latter condition means that the $k$-th row in $R$ is free.

     Set  $T(R):=T(R)+1$; this means that ``now the column $a_k=a_{k-1}$ has been used $T(R)$ times''.

  \item Otherwise, set $a_k>a_{k-1}$ to be the minimal number such that
  either (1)   there is a rook positioned at $(k,a_k)$,
 or (2) the $a_k$ column is free.

Next, we match $R$ and
     $R'=R\setminus\{(2,a_2)\}$ or $R'=R\cup\{(2,a_2)\}$ provided that  $R'$ belongs to $\Delta$.

     If $R$ is not matched, set $T(R):=1$.
\end{enumerate}

\medskip
\textbf{Remark.} If $k<n$, then $(\ast)$ guarantees that $a_k$ is well-defined. For the last row $a_n$ is ill-defined if and only if ($\ast$) is an equality and $R$ has all the rows full.

\medskip
 Eventually we have all the rows full  for non-matched simplices  (except for the unique zero-dimensional simplex).

\medskip
Now let us prove that the above defined matching is acyclic.
Take a directed path
$$R_1\nearrow Q_1 \searrow R_2 \nearrow Q_2 \searrow
\cdots . $$  Recall that
$R_i\nearrow Q_i$ if and only if $Q_i=R_i\cup \{(s_i,a_{s_i})\}$ , the
first $s_i-1$  rows of $R_i$ are full, and $a_{s_i}$ is
the first free column after $a_{s_i-1}$.

Let us prove that $(s_i,a_{s_i})$  strictly decreases along the path wrt lexicographic order.
This will imply the acyclicity.

For $Q_i \searrow R_{i+1}$, we have  $R_{i+1}=Q_{i}\setminus \{(p_i,q_i)\}$
for some $(p_i,q_i)\in Q_i$ (there are no conditions when we
remove a rook from $Q_i$). It suffices to consider the first two steps in our
directed path:
$$R_1 \nearrow Q_{1}=R_1\cup\{(s_1,a_{s_1})\}
\searrow R_2=Q_1\setminus\{(p_2,q_2)\}.$$
\begin{itemize}
    \item If $p_2>s_1$ or $p_2=s_1$ and $a_{s_1}<q_2$
    (the removed rook is below or right on $(s_1,a_{s_1})$,
    the added rook at the first
    step) our path stop, because $R_2$ is paired with
    $R_2\setminus \{(s_1,a_{s_1})\}$.
    \item  If $p_2<s_1$ or $p_2=s_1$ and $a_{s_1}>q_2$
    (the removed rook is above or left $(s_1,a_{s_1})$),
    then we have that $s_2< s_1$ or $s_2= s_1$ and
    $a_{s_2}<a_{s_1}$.
\end{itemize}

Summarizing, all critical faces (except for the unique zero-dimensional one) have all the rows full. Therefore
$\Delta^{k_1,\ldots,k_n;l_1,\ldots,l_m}_{m,n}$ is
$(k_1+k_2+\cdots+k_n-2)$-connected.  \hfill $\square$

\section{Tverberg-Van Kampen-Flores type results for $j$-wise disjoint partitions of a simplex }
\label{sec:sinisa}

Recall that a coloring of a set $S\subset \mathbb{R}^d$ is a partition $S = S_1\uplus \dots \uplus S_k$, where $S_i$ are the corresponding monochromatic sets.   By definition a subset $C\subseteq S$ is a {\em rainbow set}
if it contains at most $1$ point from each of the color classes $S_i$.
\bigskip

\begin{thm}\label{thm:application}
Let $r$ be a prime power and $j\geq 1$. Suppose that $\{S_i\}_{i=1}^k$ is a collection of $k$ finite sets of points in
$\mathbb{R}^d$ (called colors). Assume that the cardinalities $m_i = \vert S_i\vert$ satisfy the inequality
$jm_i-1\leq r$ for each $i=1,...,k$. If $(r-1)(d+1)\leq (j-1)m-1$, where $m := m_1\dots + m_k$, then it is possible to partition the set  $S = S_1\uplus \dots \uplus S_k$ into $r$ rainbow, $j$-wise disjoint sets $S = C_1\uplus\dots\uplus C_r$, so that their convex hulls
intersect,
\[
    {\rm conv}(C_1)\cap \dots \cap {\rm conv}(C_r) \neq \emptyset \, .
\]
\end{thm}

\medskip\noindent
{\bf Proof:} The rainbow sets span the multicolored simplices
which are encoded as the simplices of the simplicial complex
$([pt]_{\Delta (2)}^{*(m_1)})\ast \cdots \ast ([pt]_{\Delta
(2)}^{*(m_k)})$. Indeed these are precisely the simplices which
are allowed to have at most $1$ vertex in each of $k$ different
colors. The configuration space of all $r$-tuples of $j$-wise
disjoint multicolored simplices is the simplicial complex,

$$K=(([pt]_{\Delta (2)}^{*(m_1)})\ast \cdots \ast ([pt]_{\Delta
(2)}^{*(m_k)}))_{\Delta (j)}^{*r}$$

Since the join and deleted join commute, this complex is
isomorphic to,
$$K=([pt]_{\Delta (2)}^{*(m_1)})_{\Delta (j)}^{*r}\ast \cdots
\ast ([pt]_{\Delta (2)}^{*(m_k)})_{\Delta (j)}^{*r}$$
where $pt$ is a one-point simplicial complex.

If we suppose, contrary to the statement of the theorem, that the
intersection of images of any $r$, $j$-wise disjoint multicolored
simplices is empty, the associated mapping $F : K\rightarrow
(\mathbb{R}^d)^{*r}$ would miss the diagonal $D \subset (\mathbb{R}^d)^{*r}$. By
composing this map with the orthogonal projection to
$D^{\perp}$, and after the radial projection to the unit
sphere in $D^{\perp}$, we obtain a
$(\mathbb{Z}/p)^\alpha$-equivariant mapping,
$$\tilde F : K\rightarrow S^{(r-1)(d+1)-1}.$$
The complex $([pt]_{\Delta (2)}^{*(m_i)})_{\Delta (j)}^{*r}$ is
a multiple chessboard complex $\Delta_{m_i,r}^{1,j-1}$. Since by
assumption $jm_i-1\leq r$, this complex is $(m_i(j-1)-2)$-connected
by the main result from \cite{jvz}. Hence the complex $K$ is
$(m(j-1)-2)$-connected. By our assumption $m(j-1)-2\geq
(r-1)(d+1)-1$, so in light of Volovikov's theorem \cite{Vol96-2} such a
mapping $\tilde F$ does not exist. \hfill $\square$

\medskip
The following obvious corollary of Theorem~\ref{T:connect} is more suitable for applications in
the rest of the section.

\begin{cor}\label{dusko2}
By interchanging the rows and the columns of the multiple chessboard complex in Theorem~\ref{T:connect}, we obtain that
 the complex $\Delta_{m,n}^{k_1,...,k_n;l_1,...,l_m}$ is
$(l_1+\cdots +l_m-2)$-connected if $l_1+\cdots +l_m\leq k_1+\cdots +k_n-m+1$.
\end{cor}

\begin{thm}
Let $r$ be a prime power. Assume that positive integers $k, r, N, j$ and $d$ satisfy the inequalities $(k+1)r+r-1\leq (N+1)(j-1)$ and
$(r-1)(d+1)+1\leq r(k+1)$. Then for every continuous map $ f: \Delta^N \rightarrow \mathbb{R}^d$ there exist $r$, $j$-wise disjoint faces
of the simplex $\Delta^N$ of dimension at most $k$, whose images
have a nonempty intersection.
\end{thm}

\medskip\noindent
{\bf Proof:} The faces of dimension at most $k$ form the
$k$-skeleton $(\Delta^N)^{(k)}=[pt]_{\Delta (k+2)}^{*(N+1)}$.
The configuration space of all $r$-tuples of $j$-wise disjoint
$k$-dimensional faces of this skeleton is the simplicial complex,

$$K=([pt]_{\Delta (k+2)}^{*(N+1)})_{\Delta (j)}^{*r}.$$

This is a generalized chessboard complex
$K=\Delta_{N+1,r}^{k+1;j-1}$. Since by our assumption $(k+1)r\leq
(N+1)(j-1)-r+1$, this complex $K$ is by Corollary \ref{dusko2}
$((k+1)r-2)$-connected.

If we suppose, contrary to the statement of the theorem, that the
intersection of images of any $r$, $j$-wise disjoint
$k$-dimensional faces is empty, the associated mapping $F :
K\rightarrow (\mathbb{R}^d)^{*r}$ would miss the diagonal
$D$.

As in the proof of the previous theorem we obtain a
$(\mathbb{Z}/p)^\alpha$-equivariant mapping,

$$\tilde F : K\rightarrow S^{(r-1)(d+1)-1}.$$

We have already observed that $K$ is $((k+1)r-2)$-connected, and by our
assumption $r(k+1)-2\geq (r-1)(d+1)-1$, so in light of Volovikov's
theorem \cite{Vol96-2} such a mapping $\tilde F$ does not exist. \hfill
$\square$

\begin{thm}
Let $r$ be a prime power. Suppose that $q, r, j$ and $d$ are positive integers and let $\{S_i\}_{i=1}^k \subseteq \mathbb{R}^d$ is a
collection of colored points where all color classes $S_i$ are of the same cardinality $m$.
Then if $qr\leq m(j-1)-r+1$ and $(r-1)(d+1)+1\leq qrk$, then it is always possible to
partition the set $S := \cup_{i=1}^k S_i$  into $r$  $j$-wise disjoint sets containing
at most $q$ points of each color, so that their convex hulls ${\rm conv}(S_i)$ have a non-empty
intersection.
\end{thm}

\medskip\noindent
{\bf Proof:} The sets containing at most $q$ points of each color
span the multicolored simplices which are encoded as the simplices
of the simplicial complex $([pt]_{\Delta (q+1)}^{*m})^{*k}$.
Indeed, these are precisely the simplices which are allowed to have
at most $q$ vertices in each of $k$ different colors. The
configuration space of all $r$-tuples of $j$-wise disjoint
multicolored simplices is the simplicial complex,

$$K=(([pt]_{\Delta (q+1)}^{*m})^{*k})_{\Delta (j)}^{*r}.$$

Since the join and deleted join commute, this complex is
isomorphic to,

$$K=(([pt]_{\Delta (q+1)}^{*m})_{\Delta (j)}^{*r})^{*k}.$$

If we suppose, contrary to the statement of the theorem, that the
intersection of images of any $r$, $j$-wise disjoint multicolored
simplices is empty, the associated mapping $F : K\rightarrow
(\mathbb{R}^d)^{*r}$ would miss the diagonal $D$. As before, by
composing this map with the orthogonal projection to
$D^{\perp}$, and after the radial projection to the unit
sphere in $D^{\perp}$, we obtain a
$(\mathbb{Z}/p)^\alpha$-equivariant mapping,

$$\tilde F : K\rightarrow S^{(r-1)(d+1)-1}.$$

The complex $([pt]_{\Delta (q+1)}^{*m})_{\Delta (j)}^{*r}$ is a
multiple chessboard complex $\Delta_{m,r}^{q,j-1}$. Since we
assumed $qr\leq (j-1)m-r+1$, this complex is $(qr-2)$-connected by
Corollary \ref{dusko2}. Hence the complex $K$ is
$(qrk-2)$-connected. By our assumption $qrk\geq (r-1)(d+1)+1$, so
in light of Volovikov's theorem \cite{Vol96-2} such a mapping $\tilde F$ does
not exist. \hfill $\square$

\begin{thm}
Let $r$ be a prime power. Suppose that $q, r, j$ and $d$ are positive integers and let $\{S_i\}_{i=1}^k \subseteq \mathbb{R}^d$ is a
collection of colored points where all color classes $S_i$ are of the same cardinality $m$.
If $jm-1\leq qr$ and $(r-1)(d+1)+1\leq (j-1)mk$, then it is possible to
divide all points in $r$, $j$-wise disjoint sets containing at most
$q$ points of each color, so that their convex hulls ${\rm conv}(S_i)$  have a non-empty
intersection.
\end{thm}

\medskip\noindent
{\bf Proof:} As before the sets containing at most $q$ points of each color
span the multicolored simplices which are encoded as the simplices
of the simplicial complex $([pt]_{\Delta (q+1)}^{*m})^{*k}$.
Indeed these are precisely the simplices which are allowed to have
at most $q$ vertices in each of $k$ different colors. The
configuration space of all $r$-tuples of $j$-wise disjoint
multicolored simplices is the simplicial complex,

$$K=(([pt]_{\Delta (q+1)}^{*m})^{*k})_{\Delta (j)}^{*r}.$$

Since the join and deleted join commute, this complex is
isomorphic to,

$$K=(([pt]_{\Delta (q+1)}^{*m})_{\Delta (j)}^{*r})^{*k}.$$

If we suppose, contrary to the statement of the theorem, that the
intersection of images of any $r$, $j$-wise disjoint multicolored
simplices is empty, the associated mapping $F : K\rightarrow
(\mathbb{R}^d)^{*r}$ would miss the diagonal $D$. As before, from here by an equivariant deformation  we obtain a
$(\mathbb{Z}/p)^\alpha$-equivariant mapping,

$$\tilde F : K\rightarrow S^{(r-1)(d+1)-1} \, .$$

The complex $([pt]_{\Delta (q+1)}^{*m})_{\Delta (j)}^{*r}$ is the
multiple chessboard complex $\Delta_{m,r}^{q,j-1}$. Since we
assumed $(j-1)m\leq qr-m+1$, this complex is
$((j-1)m-2)$-connected by Corollary \ref{dusko2}. Hence the
complex $K$ is $((j-1)mk-2)$-connected. By our assumption
$(j-1)mk\geq (r-1)(d+1)+1$, and again this is in contradiction with Volovikov's theorem \cite{Vol96-2}.
 \hfill $\square$

For illustration let us consider a very special case of this theorem $q=1$ and
$j=2$.

\begin{thm}\label{thm:spec-case}
Let $r$ be a prime power. Given $k$ finite sets of points in
$\mathbb{R}^d$ (called colors), of $m$ points each, so that
$2m-1\leq r$ and $(r-1)(d+1)+1\leq mk$, it is possible to divide
the points in $r$ pairwise disjoint sets containing at most $1$
point of each color, so that their convex hulls intersect.
\end{thm}

\begin{rem}{\rm
It is easy to see that the assumptions on the total number of
points is the best possible, since the set of $(r-1)(d+1)$ points
in the general position could not be divided in $r$ disjoint sets
whose convex hulls intersect. }
\end{rem}

\subsection{A comparison with known results}

 It is interesting to compare results from the previous section with similar  results  from \cite{bfz} (Section 9). Note that the proof methods are quite different. We use high connectivity of the multiple chessboard complex, established in Section \ref{sec:dusko}, while the authors of  \cite{bfz} use the `constraint method', relying on the `optimal colored Tverberg theorem' from \cite{bmz}, as a `black box' result.

\medskip
 For illustration, let us compare our Theorem \ref{thm:spec-case} to  Theorem 9.1 from \cite{bfz}.

 \medskip Let us choose  $k\geq 2(d+1)$ in Theorem \ref{thm:spec-case} and select the smallest $m$ satisfying the inequality $(r-1)(d+1)+1\leq mk$, meaning that we are allowed to assume
 \[
 (m-1)k  <(r-1)(d+1)+1\leq mk \, .
 \]
 From here we immediately deduce the inequality $2m-1\leq r$ and, as a consequence of Theorem \ref{thm:spec-case}, we have the following result.

 \begin{cor}\label{cor:cor}
   Let $r$ be a prime power.  Assume $k\geq 2(d+1)$  and choose $m$ satisfying the inequality $(r-1)(d+1)+1\leq mk$. Suppose that $S\subset \mathbb{R}^d$ is a set of cardinality $mk$,  evenly colored by $k$ colors (meaning that  $S = \cup_{i=1}^k~S_i $  where $\vert S_i\vert = m$ for each $i$). Then   it is possible to select
 $r$ pairwise disjoint subsets $C_i\subset S$, containing at most $1$
point of each color, so that $\cap_{i=1}^r {\rm conv}(C_i)\neq\emptyset$.
 \end{cor}

 This result clearly follows from  Theorem 9.1 if we assume that $r$ is a prime.  Corollary \ref{cor:cor} illustrates the phenomenon that there exist instances of the `optimal colored Tverberg theorem' (Theorem 9.1 in \cite{bfz}) which remain valid if the condition on $r$ being a prime is relaxed to $r$ is a prime power.

\subsection{A remark on Tverberg A-P conjecture}

In this section we briefly discuss the problem whether each admissible
r-tuple is Tverberg prescribable. This problem, as formulated in
\cite{bfz}, will be referred to as the Tverberg A-P problem or the
Tverberg A-P conjecture.

\begin{defin}
For $d\geq 1$ and $r\geq 2$, an $r$-tuple $d = (d_1,...,d_r)$ of
integers is admissible if, $[\frac d2]\leq d_i\leq d$ for all $i$,
and $\sum_{i=1}^r (d -d_i)\leq d$. An admissible r-tuple is
Tverberg prescribable if there is an $N$ such that for every
continuous map $f : \Delta^N\rightarrow {\mathbb R}^d$ there is a
Tverberg partition $\{\sigma_1,...,\sigma_r\}$ for $f$ with $\dim
(\sigma_i) = d_i$.
\end{defin}

\noindent {\bf Question:} (Tverberg A-P problem; \cite{bfz}
(Question 6.9.)) Is every admissible $r$-tuple Tverberg
prescribable?

\medskip
As shown in \cite{F}, (Theorem 2.8.), the answer to the above
question is negative. It was also demonstrated that a more realistic conjecture  arises if the condition $[\frac d2]\leq d_i\leq d$, in
the definition of admissible $r$-tuple, is replaced by a stronger requirement
$\frac {(r-1)}r(d-1)\leq d_i\leq d$ for all $i$.

\medskip

Here we remark that a positive answer to the modified question is quite straightforward in the case $r\geq d$.
Indeed, in this case we have for all $i$
$$d_i\geq \frac {(r-1)}r(d-1)\geq d-1-\frac {(d-1)}r>d-2.$$
So, in this case each $d_i$ is equal to either $d-1$ or $d$, and the A-P conjecture reduces to the `balanced case', established in \cite{jvz2}.

\end{document}